\documentclass{amsart}
\begin{document}

\vfuzz2pt 
\hfuzz2pt 
\newtheorem{thm}{Theorem}[section]
\newtheorem{corollary}[thm]{Corollary}
\newtheorem{lemma}[thm]{Lemma}
\newtheorem{proposition}[thm]{Proposition}
\newtheorem{defn}[thm]{Definition}
\newtheorem{remark}[thm]{Remark}
\newtheorem{example}[thm]{Example}
\newtheorem{fact}[thm]{Fact}
\
\newcommand{\norm}[1]{\left\Vert#1\right\Vert}
\newcommand{\abs}[1]{\left\vert#1\right\vert}
\newcommand{\set}[1]{\left\{#1\right\}}
\newcommand{\Real}{\mathbb R}
\newcommand{\eps}{\varepsilon}
\newcommand{\To}{\longrightarrow}
\newcommand{\BX}{\mathbf{B}(X)}
\newcommand{\A}{\mathcal{A}}

\def\proof{\medskip Proof.\ }
\font\lasek=lasy10 \chardef\kwadrat="32 
\def\kwadracik{{\lasek\kwadrat}}
\def\koniec{\hfill\lower 2pt\hbox{\kwadracik}\medskip}

\newcommand*{\C}{\mathbf{C}}
\newcommand*{\R}{\mathbf{R}}
\newcommand*{\Z}{\mathbf {Z}}

\def\sb{f:M\longrightarrow \C ^n}
\def\det{\hbox{\rm det}\, }
\def\detc{\hbox{\rm det }_{\C}}
\def\i{\hbox{\rm i}}
\def\tr{\hbox{\rm tr}\, }
\def\rk{\hbox{\rm rk}\,}
\def\vol{\hbox{\rm vol}\,}
\def\Im {\hbox{\rm Im}\, }
\def\Re{\hbox{\rm Re}\, }
\def\interior{\hbox{\rm int}\, }
\def\e{\hbox{\rm e}}
\def\pu{\partial _u}
\def\pv{\partial _v}
\def\pui{\partial _{u_i}}
\def\puj{\partial _{u_j}}
\def\puk{\partial {u_k}}
\def\div{\hbox{\rm div}\,}
\def\Ric{\hbox{\rm Ric}\,}
\def\r#1{(\ref{#1})}
\def\ker{\hbox{\rm ker}\,}
\def\im{\hbox{\rm im}\, }
\def\I{\hbox{\rm I}\,}
\def\id{\hbox{\rm id}\,}
\def\exp{\hbox{{\rm exp}^{\tilde\nabla}}\.}
\def\cka{{\mathcal C}^{k,a}}
\def\ckplusja{{\mathcal C}^{k+1,a}}
\def\cja{{\mathcal C}^{1,a}}
\def\cda{{\mathcal C}^{2,a}}
\def\cta{{\mathcal C}^{3,a}}
\def\c0a{{\mathcal C}^{0,a}}
\def\f0{{\mathcal F}^{0}}
\def\fnj{{\mathcal F}^{n-1}}
\def\fn{{\mathcal F}^{n}}
\def\fnd{{\mathcal F}^{n-2}}
\def\Hn{{\mathcal H}^n}
\def\Hnj{{\mathcal H}^{n-1}}
\def\emb{\mathcal C^{\infty}_{emb}(M,N)}
\def\M{\mathcal M}
\def\Ef{\mathcal E _f}
\def\Eg{\mathcal E _g}
\def\Nf{\mathcal N _f}
\def\Ng{\mathcal N _g}
\def\Tf{\mathcal T _f}
\def\Tg{\mathcal T _g}
\def\diff{{\mathcal Diff}^{\infty}(M)}
\def\embM{\mathcal C^{\infty}_{emb}(M,M)}
\def\U1f{{\mathcal U}^1 _f}
\def\Uf{{\mathcal U} _f}
\def\Ug{{\mathcal U} _g}
\def\[f]{{\mathcal U}^1 _{[f]}}

\title{A moduli space of minimal affine Lagrangian submanifolds}
\author{Barbara Opozda}

\subjclass{ Primary: 53C40, 57R40, 58B99 Secondary: 53C38, 58A10}

\keywords{infinite dimensional Fr\'echet manifold, variation,  phase
function, tubular mapping, H\"older-Banach space, differential
operator, elliptic regularity theorems}

 \thanks{The research partially supported by the
grant NN 201 545738} \maketitle

\address{Instytut Matematyki UJ, ul. \L ojasiewicza  6, 30-348 Cracow,
Poland}

\email{Barbara.Opozda@im.uj.edu.pl}

\vskip 1in \noindent

\vskip 0.5in \noindent
 {\bf Abstract.}
 It is proved that the moduli space of all connected compact orientable embedded minimal
 affine Lagrangian submanifolds of a complex equiaffine space constitutes an infinite dimensional
 Fr\'echet manifold (if it is not $\emptyset$). The moduli space of
 all connected compact orientable metric Lagrangian embedded  surfaces in
 an almost K\"ahler 4-dimensional manifold forms an infinite dimensional
 Fr\'echet manifold (if it is not $\emptyset$).

\maketitle

\bigskip

\section{Introduction}

R. McLean  proved in \cite{McL} that special Lagrangian submanifolds
near a compact special Lagrangian submanifold of a Calabi-Yau
manifold form a  manifold of dimension $b_1$, where $b_1$ is the
first Betti number of the submanifold. Then  a few papers giving
generalizations to the cases where the ambient space is not a
Calabi-Yau  manifold but a more general type of space have been
published. All those cases are, in fact,  within metric geometry.
The aim of this paper is to prove a similar result  in the
non-metric case. Moreover, we prove a global result, that is, we
describe the set of all minimal affine Lagrangian embeddings of a
compact manifold. It turns out that this set has a nice structure.
Namely, it is an infinite dimensional Fr\'echet manifold  modeled on
the Fr\'echet space of all closed $(n-1)$-forms on the submanifold,
where $n$ is the complex dimension of the ambient space.
The main result of this paper says  that the set of all minimal
affine Lagrangian embeddings of a compact manifold into an
equiaffine complex space is a submanifold of the Fr\'echet manifold
of all compact
 submanifolds of the  complex equiaffine space. We provide a rigorous proof of this fact.

 It seems that from a differential geometry viewpoint
non-metric analogues of Calabi-Yau manifolds are equiaffine  complex
manifolds, that is, complex manifolds equipped with a torsion-free
complex connection and a non-vanishing covariant constant complex
volume form. There are many very natural complex equiaffine
manifolds. For instance, complex affine hyperspheres of the complex
affine space $\C ^n$ with an induced equiaffine structure obtained
in a way standard in affine differential geometry (see \cite{DV})
are examples. Equiaffine structures are, in general, non-metrizable.
For instance, the complex hyperspheres of $\C ^n$ with  the induced
equiaffine structure are non-metrizable. In particular, they are not
related to Stentzel's metric.

 If $N$ is a complex $n$-dimensional  space with  a complex
structure $J$, then  an $n$-dimensional real submanifold $M$ of $N$
is affine Lagrangian if $JTM$ is transversal to $M$.  Of course, if
$N$ is  almost Hermitian, then Lagrangian (in the metric sense)
submanifolds, for which $JTM$ is orthogonal to $TM$, are affine Lagrangian and there are many affine
Lagrangian submanifolds which are not metric Lagrangian even if
the ambient space is almost Hermitian.

In order to discuss  minimality of  submanifolds  a metric structure
is not necessary.  It is sufficient to have  induced volume elements
on submanifolds. Such a situation  exists in   case of affine
Lagrangian submanifolds.  In this case there does not exist (in
general) any mean curvature vector but there exists the Maslov
1-form which  can play, in some situations, a role similar to that
played by the mean curvature vector. Note that in the general affine
case we do not have  any  canonical duality between tangent vectors
and 1-forms. The vanishing of the Maslov form  implies that the
submanifold  is a   point where a naturally defined volume
functional attains its minimum for compactly supported variations.
Affine Lagrangian submanifolds have a phase function. It turns out
that a connected affine Lagrangian submanifold is minimal if and
only if its phase function  is constant.
If a connected affine Lagrangian submanifold is minimal (i.e. of
constant phase), then after rescaling the complex volume form in the
ambient space we can assume that the constant phase function
vanishes on $M$. Analogously to the metric case an affine Lagrangian
submanifold is called special if its phase function vanishes on $M$.
The notion of special submanifolds corresponds to the notion of
calibrations. Calibrations in Riemannian geometry were introduced in
the famous paper \cite{HL}. The notion can be generalized to the
affine case and, like in the metric case, an affine Lagrangian
submanifold is special if and only if it is calibrated by the real
part of the complex volume form in the ambient space. The minimality
of affine Lagrangian submanifolds is discussed in \cite{O_1} and
\cite{O}.

In this paper we try to assume as little as possible. In particular,
we do not assume that the ambient space $N$ is complex equiaffine
 but
we only assume that it is almost complex and endowed with a nowhere
vanishing closed complex $n$--form $\Omega$, where $2n=\dim _{\R}N$.
Then affine Lagrangian immersions $f:M\to N$, where $\dim M=n$, are
those for which $f^*\Omega \ne 0$ at each point of $M$. If $M$ is
oriented, then $\Omega$ induces on $M$ a unique volume element $\nu
$. We have $f^*\Omega=\e ^{\i\theta}\nu $, where $\theta $ is the
phase function of $f$. In this paper minimal (relative to $\Omega$)
affine Lagrangian submanifolds will be those (by definition) which
have constant phase.

We shall prove the following theorem.

\begin{thm}\label{main} Let $M$ be a connected compact oriented $n$-dimensional real
manifold admitting a minimal affine Lagrangian embedding into an
almost complex $2n$--dimensional manifold $N$ equipped with a
nowhere-vanishing complex closed $n$-form. The set of all minimal
affine Lagrangian embeddings of $M$ into $N$ has a structure of an
infinite dimensional manifold modeled on the Fr\'echet vector space
$\mathcal C ^{\infty}(\mathcal F_{closed}^{n-1})$ of all smooth
closed $(n-1)$--forms on $M$.
\end{thm}
 A  precise formulation of this theorem is given in Section 4.

 The Fr\'echet manifold in Theorem \ref{main} may have many connected
 components.
In the above theorem a manifold $M$ is fixed. But the theorem says,
in fact, about all compact (connected oriented) minimal affine
Lagrangian embedded submanifolds. Non-diffeomorphic submanifolds are
in different connected components.


At the end of this paper we observe that, in contrast with the
metric geometry, in the affine case there exist non-smooth minimal
submanifolds of smooth manifolds.
Almost the same consideration as in the proof of Theorem \ref{main}
gives the statement saying that the set of all (metric) Lagrangian
embeddings of a connected compact 2-dimensional manifold into a
4-dimensional almost K\"ahler manifold forms an infinite dimensional
Fr\'echet manifold modeled on the Fr\"echet vector space $\mathcal C
^{\infty} (\mathcal F ^1_{closed})$ (if it is not $\emptyset$).

We also give a  simple application of Theorem \ref{main} in the case
where the ambient space is the tangent bundle of a flat manifold.

\bigskip

\section{Basic notions}

Let $N$ be a $2n$--dimensional almost complex manifold with an
almost complex structure $J$.  Let $M$ be a connected
$n$-dimensional manifold and $f:M\to N$ be an immersion. We say that
$f$ is affine Lagrangian (some authors call it purely real or
totally real) if the bundle $Jf_*(TM)$ is transversal to $f_*(TM)$.
We shall call this transversal bundle the normal bundle.
The almost complex structure $J$ gives  an isomorphism between the normal bundle and the tangent bundle $TM$.
 If $\Omega$ is a nowhere vanishing complex
$n$-form on $N$, then $f$ is affine Lagrangian if and only if
$f^*\Omega \ne 0$ at each point of $M$. If $f:M \to N$ is such a
mapping that $f^*\Omega \ne 0$ at each point of $M$, then $f$ is
automatically an immersion.

Recall now the notion of a phase. Let $\bf V$ be an $n$-dimensional
complex vector space with a complex volume form $\Omega$ and $\bf U$
be its $n$-dimensional   real oriented vector subspace such that
$\Omega _{| {\bf U}}\ne 0$. Let $X_1,...., X_n$ be a positively
oriented basis of $\bf U$. Then $\Omega (X_1,...,X_n)=\mu
e^{\i\theta}$, where $\mu \in \R^+$ and $\theta \in \R$. If we
change the  basis $X_1,..., X_n$  to another positively oriented
basis of $\bf U$, then $e^{\i\theta}$ remains unchanged. $\theta$ is
called the phase or the  angle of the subspace  $\bf U$.

Assume $N$ is endowed with a nowhere vanishing complex volume form $\Omega$ and $M$ is oriented.
 For an
affine Lagrangian immersion $f:M\longrightarrow N$, at each point
$x$ of $M$ we have the phase $\theta _x$ of the tangent vector
subspace $f_*(T_xM)$ of $T_{f(x)} N$. The phase function
$x\longrightarrow \theta _x$ is multi-valued. In general, if we want
to have the phase function  to be a smooth function, it is defined
only locally.  For each point $x\in M$ there is a smooth phase
function of $f$ defined around $x$.
 The constancy of the
phase function is a well defined global notion, that is, if $\theta$
is locally constant, then it can be chosen
globally constant.

 Recall few facts  concerning the situation where the ambient space is complex equiaffine
or, in the metric case, Calabi-Yau.  A Lagrangian submanifold
(affine or metric) is minimal if and only if  it is volume
minimizing for compactly supported variations. This is equivalent to
the fact that the Maslov form vanishes. Moreover, Lagrangian
submanifolds (affine or metric) are minimal if and only if they have
constant phase.

In this paper, where, in general,  we do not assume that the ambient space is complex equiaffine, we shall say (by
definition) that an affine Lagrangian submanifold $f$ is minimal  if
and only if its phase is constant. As usual, if the phase constantly  vanishes
on $M$, then the submanifold will be  called special. If the phase $\theta$ is
constant, then we can rescale $\Omega$ in the ambient space by
multiplying it by $\e ^{-\i \theta}$ and after this change the given
immersion becomes special. But if we have a family of minimal affine
Lagrangian immersions of $M$ into $N$,  and we adjust the complex
volume form $\Omega$ to one member of the family then, in general,
the rest of the family remain only minimal.

For an oriented affine Lagrangian immersion $f:M\to  N$ we have the
induced volume form $\nu $  on $M$ defined by the condition $\nu
(X_1,..., X_n)=\vert\omega (X_1,..., X_n)\vert$, where $X_1,...,X_n$
is a positively oriented basis of $T_xM$, $x\in M$ and $\omega
=f^*\Omega$. The form $\omega$ is a real complex--valued $n$--form
on $M$. We have
$$\omega =\e ^{\i\theta } \nu ,$$
where $\theta$ is the phase function. Note that by multiplying
$\Omega$ by $\e ^{\i \alpha}$ for any $\alpha \in\R$ we do not
change the induced volume form on $M$. Decompose $\omega$ into the
real and imaginary parts:  $\omega= \omega _1+\i\omega_2$, where
$\omega _1=\cos \theta\,\nu $, $\omega _2=\sin\theta\,\nu.$ If
$W\in\mathcal X (M)$, then, since $\Omega$ is complex, we have
$$
f^*(\iota _{(Jf_*W)} \Omega)=-\iota _W\omega _2 +\i \iota _W\omega
_1,
$$
where $\iota$ stands for the interior product operator.
Hence, if  $f$ is special (i.e. $\nu =\omega_1$) we get
\begin{equation}\label{przedGriffiths}
f^*(\iota _{(J f_*W )}\Im\Omega) = \iota _W\nu.
\end{equation}

Assume now that $f_t$, $|t|<\varepsilon$, is a smooth  variation of $f$. Denote by
$\mathcal V (t, x)$ its variation vector field. Assume it is normal
to $f$ at $t=0$. Then $V:=\mathcal V_{\vert{\{0\}}\times M}$ is
equal to $Jf_*W$ for some $W\in\mathcal X (M)$. If $f$ is special and $\Omega$ is closed,
then using formula
(\ref{przedGriffiths}) and Proposition (I.b.5) from \cite{G}, we
obtain

\begin{equation}\label{Griffiths}
{d\over{dt}}\left (f^*_t\Im \Omega\right )_{|t=0} =d(\iota _W\nu).
\end{equation}
This formula is also
directly computed in \cite{O}, but there the form $\Omega$ is
assumed to be parallel relative to a torsion-free complex
connection.

We shall now give a justification of the term ``minimal'' adopted in
this paper. Assume $M$ is compact. If  affine Lagrangian immersions
$f, \tilde f: M\to N$ are cohomologous (in particular, if they are
homotopic), then the cohomology class of $\omega _i$ is equal to the
cohomology class of $\tilde\omega _i$, for $i=1,2$, where
$\tilde\omega _1=\cos \tilde\theta\, \tilde\nu$, $\tilde\omega
_2=\sin\tilde\theta\, \tilde \nu$ are the real and imaginary parts
of $\tilde{\omega} =\tilde f ^*\Omega$ and $\tilde\theta$,
$\tilde\nu$ are the phase and the induced volume element for $\tilde
f$.
 Assume that $f$ is special.
Then  $\omega _1=\nu $ and consequently
$$ \int _M \nu=\int _M \omega _1 =\int _M \tilde\omega _1=\int _M
\cos\tilde\theta \, \tilde \nu \le\int _M\tilde\nu,$$
 which means that
with the definition of minimality we adopted in this paper  compact
special (and consequently minimal) affine Lagrangian submanifolds
are volume minimizing in their respective  cohomology classes.

 Assume additionally
that $\tilde f$ is minimal with the constant phase $\tilde\theta$.
We have
$$0=\int _M\omega _2=\int _M \tilde\omega_2= \int _M \sin\tilde\theta\,
\tilde \nu =\sin\tilde\theta \int_M \tilde\nu ,$$ which means that
$\tilde\omega _2=0$, that is, $\tilde f$ is also special.

If $f$ is minimal  (special), then for any diffeomorphism $\varphi$
of $M$ $f\circ\varphi$ is minimal (special).
\bigskip

\section{Moduli spaces of compact embedded submanifolds}
 Assume first that $M$  and $N$ are arbitrary manifolds
such that $\dim M\le \dim N$. Assume  moreover that $M$ is connected
compact and it admits an embedding into $N$. Denote by $\mathcal
C^{\infty}_{emb}(M,N)$ the set of all embeddings from $M$ into $N$.
This is a well known topological space forming an open subset (in
the $\mathcal C^1$ topology) of $\mathcal C ^{\infty}(M,N)$.

 Denote by
$\mathcal M$ the space $\emb _{/ Diff ^{\infty}(M)}$ with the
quotient  topology. The equivalence class of $f\in \emb$ will be
denoted by $[f]$. For $f,g\in \emb$ we have that $f\sim g$ if and
only if the images of $f$ and $g$ are equal in $N$.

We shall now introduce a structure of an infinite dimensional
manifold (modeled on Fr\'echet spaces) on $\M$. It is certainly well
known but we have not found suitable references and moreover we need
the construction. We use the notion of  a  manifold modeled on
Fr\'echet vector spaces given in \cite{H}. We denote the  Fr\'echet
space of all  $\mathcal C ^{\infty}$ sections of a vector bundle
$E\to M$ by $\mathcal C^{\infty}(M\leftarrow E)$. Analogously the
Banach spaces of all $\mathcal C ^k$ sections of a vector bundle
$E\to M$ will be denoted by $\mathcal C ^k(M\leftarrow E)$.

The basic tool in the construction are tubular mappings. We use the
following setting of this notion. Assume that $\mathcal N _f$ is any
smooth transversal bundle for an embedding $f:M\to N$. Having any
connection on $N$ we have the exponential mapping $exp $ given by
the connection. No relation between the connection and the
transversal bundle is needed. From the theory of connections one
knows that there is an open neighborhood $\mathcal U$ of the
zero-section in the total space $\mathcal N _f$ and an open
neighborhood $\mathcal T$ of $f(M)$ in $N$ such that $exp_{\vert
\mathcal U} :\mathcal U\to \mathcal T$ is a diffeomorphism, $exp
_{|M}=\id _M$ and the differential $exp _*: T_{0_x}(\mathcal
N_f)=f_*(T_xM)\oplus \mathcal (N_f)_x\to T_xN$ of $exp$ at $0$ is
the identity for each point $x$ of $M$. The mapping $exp_{\vert
\mathcal U}$ is a tubular mapping. In order to reduce a play with
neighborhoods we shall use the following lemma, which allows to have
the whole total space $\Nf$ as the domain of a tubular mapping. In
what follows $\mathcal N _f$ will denote either the transversal
vector bundle or its total space depending on the context.

\begin{lemma}\label{lemacik_o_funkcji}
Let $E\longrightarrow M$ be a Riemannian vector bundle and $\mathcal
U_{\varepsilon}$ be the neighbourhood of the zero section of $E$
given as follows
$$\mathcal U_{\varepsilon}=\{v\in E;\ \mid v\mid<\varepsilon\}, $$
where $\mid\ \mid$ is the norm on fibers of $E$ determined by the
Riemannian structure. There is a fiber-respecting diffeomorphism
$\sigma :E\longrightarrow \mathcal U_{\varepsilon}$ which is  the
identity on $\mathcal U_{\varepsilon/2}$.
\end{lemma}

\proof Let $\psi : [0, \infty)\to \R$ be a smooth function such that
$\psi (t)=t$ for $t\le \varepsilon/2$, $\psi (t)\le \varepsilon$ for
$t>\varepsilon/2$ and $\psi(t)\rightarrow \varepsilon$ for
$t\to{\infty}$. Then the function $\Upsilon(t)=(1/t)\psi (t)$ is
also a smooth function on $[0,\infty)$. The mapping
$\sigma:E\to\mathcal U_{\varepsilon} $ given by
$$\sigma (v)=\Upsilon (\mid v\mid )v$$
satisfies the required conditions.  \koniec

We now endow the bundle $\mathcal N _f$ with any Riemannian metric.
Since $M$ is compact, there is $\varepsilon>0$ such that $\mathcal U
_{\varepsilon}\subset \mathcal U$. We use Lemma
\ref{lemacik_o_funkcji} for $\mathcal U_{\varepsilon}$ and take the
tubular mapping $\Ef = exp \circ \sigma $.  The tubular neighborhood
$\Ef(\Nf)$ of $f(M)$ will be denoted by $\Tf$.

The set $\mathcal C^{\infty}_{emb}(M,\mathcal T_f )$ is open in
$\emb$ (in the $\mathcal C ^0$-topology).
 Consider the mapping
\begin{equation}
\Psi: \mathcal C^{\infty}_{emb}(M,\mathcal T_f )\ni h\longrightarrow
\Pi_f\circ\Ef ^{-1}\circ h\in  \mathcal C^{\infty}(M,M),
\end{equation}
where $\Pi _f:\Nf\longrightarrow M$ is the  natural projection. The
mapping is continuous and the set $\diff $ is open in $\mathcal C
^{\infty} (M,M)$ (in the $\mathcal C ^1$-topology). Thus the set
\begin{equation}
{\mathcal U}^1_f =\Psi^{-1}(\diff )= \{h\in\mathcal C^{\infty}(M,\Tf
); \ \Pi _f\circ \Ef ^{-1}\circ h\in\diff \}
\end{equation}
is open in $\emb$ in the $\mathcal C^1$ topology. Observe that
$h\in\mathcal U ^1_f$ if and only if there is a section $V\in
{\mathcal C}^{\infty}(M\longleftarrow\Nf )$ and $\varphi \in \diff $
such that
\begin{equation}
\Ef\circ V=h\circ \varphi .
\end{equation}
The set $\U1f$ has the following properties:

1) If $h\in\U1f$ and $\varphi \in\diff$, then $h\circ\varphi
\in\U1f$ .

2) For every $\varphi \in\diff$ we have ${\mathcal U}
^1_{f\circ\varphi}=\U1f$.

Take the neighborhood $\[f]= \{[h]\in\M ; \ h\in\U1f\}$ of $[f]$ in
$\M$. Observe that the elements of $\[f]$ can be parametrized
simultaneously. Namely, we have
\begin{lemma}\label{simultaneously}
Let $\xi _0 \in \[f] $ and $h_0\in\emb$ be its fixed
parametrization. For each $\xi \in \[f]$ there is a unique
parametrization $h_{\xi}\in\emb$ of $\xi$ such that
$$\Pi _f\circ\Ef ^{-1}\circ h_0 =\Pi _f\circ \Ef ^{-1}\circ
h_{\xi}$$
\end{lemma}
\proof We first reparametrize $f$ in such a way that after the
reparametrization
$$\Pi _f \circ\Ef ^{-1}\circ h_0=\id _M.$$
Assume that $f$  is  already parametrized in this way.
 For every $h\in \U1f $ the mapping $\varphi ^{-1}=\Pi _f\circ\Ef ^{-1}\circ
 h$ is a diffeomorphism and it is sufficient to replace $h$ representing $[h]$ by
 $h\circ\varphi$. The uniqueness is obvious.
 \kwadracik

By the above lemma we see that $\[f]$ can be identified with the set
\begin{equation}{\mathcal U_{[f]}} =\{h\in\mathcal C^{\infty}_{emb}(M,\mathcal T_f ) ;
 \ \Pi _f\circ\Ef ^{-1}\circ h=\id_M\}.\end{equation}
We now define the bijection
$$
u_{[f]}:\mathcal U_{[f]}\longrightarrow \C ^{\infty}(M\leftarrow
\Nf)$$ as follows:
\begin{equation}
u_{[f]}(h)\longrightarrow \Ef^{-1}\circ h .
\end{equation}
We see that
$$u_{[f]} ^{-1} (V)=\Ef\circ V$$
and $\Ef\circ V$ has values in $\Tf$. If $ U$ is an open subset of
$\Tf$, then
$$u_{[f]}(\{h\in \mathcal U_{[f]};  h(M)\subset
U\})=\{V\in\mathcal C ^{\infty}(M\leftarrow\Nf);\
V(M)\subset\Ef^{-1}( U)\}$$ and hence  is open in $\mathcal C
^{\infty}(M\leftarrow \Nf)$.

 Assume now that $f, g\in \emb$ and $\mathcal
U_{[f]}\cap \mathcal U_{[g]}\ne\emptyset$. Take $\xi _0\in \mathcal
U_{[f]}\cap \mathcal U_{[g]}$ and  fix its parametrization $h_0$.
Reparametrize $f$ and $g$ as in Lemma \ref{simultaneously} adjusting
the parametrizations to $h_0$. Then
\begin{eqnarray*}\mathcal U_{[f]}\cap \mathcal U_{[g]}&=\{h\in \mathcal C
^{\infty}_{emb}(M,\Tf\cap\Tg);\ \Pi _f\circ \Ef ^{-1}\circ h=\id _M,
\ \Pi _g\circ\Eg ^{-1}\circ h=\id _M\}\\
&=\{\Ef\circ V;\ V\in   \mathcal C^{\infty} (M\leftarrow\Nf);\ V(M)\subset \Ef ^{-1} (\Tf\cap\Tg )\}\\
\end{eqnarray*}
and consequently
$$u_{[f]}(\mathcal U_{[f]}\cap \mathcal U _{[g]})=\{
\ V\in \mathcal C^{\infty} (M\leftarrow\Nf );\  V(M)\subset \Ef
^{-1} (\Tf\cap\Tg )\}$$

 The mapping $\Eg^{-1}\circ\Ef
 :\Ef^{-1}(\Tf\cap\Tg)\to \Eg ^{-1}(\Tf\cap\Tg )$ is  smooth and
 fiber respecting  (because of specially chosen parametrizations  $f$ and $g$).
 It is known,  \cite{H},
  that the set $u_{[f]}(\mathcal U_{[f]}\cap \mathcal U _{[g]})$ is open in the Fr\'echet space $\mathcal C^{\infty}(M\leftarrow \mathcal N _f)$
and the mapping
\begin{eqnarray*}
&u_{[g][f]}:u_{[f]}(\mathcal U_{[f]}\cap \mathcal U _{[g]})\ni V\to
\Eg ^{-1}\circ\Ef\circ V \in u_{[g]}(\mathcal U_{[f]}\cap \mathcal U
_{[g]})
\end{eqnarray*}
is  smooth. For the same reason the  set  $u_{[g]}(\mathcal
U_{[f]}\cap \mathcal U_{[f]})$ is open and  the mapping $u_{[f][g]}$
 is  smooth.

We have
built a smooth atlas on $\M$. Hence we have
\begin{thm}
Let $M$ be a connected compact manifold admitting an embedding in  a
manifold $N$. Then  $\mathcal M$ is an infinite  dimensional
manifold modeled on the Fr\'echet  vector spaces $\mathcal C
^{\infty}(M\leftarrow \Nf)$ for $f\in\mathcal C_{emb}
^{\infty}(M,N)$, where  $\mathcal N _f$ is any bundle transversal to
$f$.
\end{thm}
In the theorem $\Nf$ can be replaced by  any bundle isomorphic (over
the identity on $M$) to the transversal bundle $\Nf$.

In what follows the Fr\'echet space of all smooth $r$-forms on $M$
will be denoted by $\mathcal C ^{\infty}(\mathcal F ^{r})$. The
Banach space of $r$-forms   of class $\mathcal C ^k$,
$k\in\mathbf N$, will be denoted by $\mathcal C ^k(\mathcal F ^r)$.

Assume now additionally that $N$ is a $2n$-dimensional manifold with
an almost complex structure $J$ and $M$ is $n$-dimensional
orientable with a fixed volume form $\nu$. Having the volume element
$\nu$, we have an isomorphic correspondence between tangent vectors
and $(n-1)$-forms. It is given by the interior multiplication
$$T_xM\ni W\longrightarrow \iota _W\nu \in  \Lambda ^{n-1}(T_xM)^*,$$
for $x\in M$. If $f:M\to N$ is affine Lagrangian, then by composing
this isomorphism with the isomorphism determined by $J$ between  the
tangent bundle $TM$ and the normal bundle $\mathcal N _f$ we get an
isomorphism, say $\rho$, of vector bundles
\begin{equation}\label{rho}
\rho: \Lambda ^{n-1} TM^*\longrightarrow \mathcal N _f.
\end{equation}
The isomorphism gives a  smooth isomorphism (linear smooth
diffeomorphism) $\wp $ between Fr\'echet vector spaces $\mathcal C
^{\infty}(\mathcal F ^{n-1})$ and $\mathcal C^ {\infty}(M\leftarrow
\Nf)$  given by $\wp (\gamma)=\rho\circ \gamma$.

We now have
\begin{thm}\label{Mal}
Let $M$ be a connected compact orientable $n$-dimensional real
manifold admitting an affine Lagrangian embedding into a
$2n$-dimensional almost complex manifold $N$. The set $\M aL
=\{[f]\in \M;\ f \ is \ affine\ Lagrangian\}$ is an infinite
dimensional manifold modeled on the Fr\'echet vector space
$\C^{\infty}(\mathcal F^{n-1})$.
\end{thm}
\proof
For each $y\in N$ there is an open neighborhood $U_y$ of $y$  in $N$ and
a smooth complex $n$-form $\Omega _y$ on $N$ such that $\Omega _y\ne 0$ at each point of $U_y$.
Let $U_{y_1},..., U_{y_l}$ cover $f(M)$. Set $\tilde{\Theta} _j=  {\mathcal E _f} ^*\Omega _{y_j} $.
 Consider the mapping
\begin{equation}\mathcal C ^1(M\leftarrow \Nf)\ni V\to (V^*\tilde{\Theta}_1,...,V^* \tilde{\Theta}_l)\in \mathcal
(C ^0(\mathcal F (\C)))^l.
\end{equation}
where  $\mathcal C ^0(\mathcal F(\C))$ stands for the space of all
real complex-valued $n$-forms on $M$ of class $\mathcal C ^0$. It is
known, see Theorem 2.2.15 from \cite{Ba}, that  this mapping is
continuous between Banach spaces. Hence
$$\tilde{\mathcal U}=\{V\in\mathcal C ^{\infty}(M\leftarrow\Nf);\
((V^*\tilde{\Theta }_1)_x, ... , (V^*\tilde {\Theta }_l)_x)\ne 0 \
\forall x\in M\}$$ is open in $\mathcal C^{\infty}(M\leftarrow\Nf)$.
It is clear that $[h]\in \mathcal U ^1_{[f]}$ is affine Lagrangian
if and only $u_{[f]}([h])\in \tilde {\mathcal U}$. We now compose
$u_{[f]}$ with the  isomorphism $\wp ^{-1}$, where $\wp$ is
determined by any fixed volume form on $M$. \koniec

In the above atlas we can compose a chart $u_{[f]}$ with a bijective
mapping, say $\phi$, sending an open neighborhood of $0$ in
$\C^{\infty}(\mathcal F^{n-1})$ onto an open neighborhood of $0$ in
$\C^{\infty}(\mathcal F^{n-1})$ and such that $\phi $ and $\phi
^{-1}$ are smooth in the sense of the theory of Fr\'echet vector
spaces. This does not change the differentiable structure on
$\mathcal M aL$. We shall use this possibility in the next section.

\bigskip
\section{The moduli space of minimal submanifolds}

The precise formulation of Theorem \ref{main} is the following

\begin{thm}\label{main-theorem} Let $N$ be a $2n$-dimensional  almost complex manifold equipped with a
smooth nowhere-vanishing closed complex $n$-form $\Omega$. Let $M$
be a connected compact oriented    $n$-dimensional real manifold
admitting a minimal (relative to $\Omega$)  affine Lagrangian
embedding into  $N$.
 Then the set $$\M maL= \{[f]\in \M aL;\  f \ is \
minimal\}$$ is an infinite dimensional manifold modeled on the
Fr\'echet  vector space \newline $\C^{\infty}(\mathcal
F^{n-1}_{closed})$. It is a submanifold of $\M aL$.
\end{thm}



\proof
 We shall improve the charts obtained in Theorem
\ref{Mal} in such a way that the set $\mathcal M maL$ will get a
structure of a submanifold of $\mathcal MaL$ in the sense of the theory of Fr\'echet manifolds.
 Let $f:M\to N$ be a
given minimal affine Lagrangian embedding. By rescaling $\Omega$ in
the ambient space we  make $f$  special. We have  the normal bundle
$\Nf=J f_*(TM)$. Fix a tubular mapping  $\Ef :\Nf\to\Tf$. For each
section $V\in \mathcal C ^{\infty}(M\leftarrow \Nf)$ we have the
embedding $f_V=\Ef\circ V$. In general, $f_V$ is neither special nor
minimal nor even affine Lagrangian.  Consider the mapping
$$\tilde P: \mathcal C ^{\infty}(M\leftarrow \Nf)\ni V\to \tilde P (V)=f_V^* (\Im\,
\Omega)\in\mathcal C ^{\infty} (\mathcal F^{n}).$$ Of course $\tilde
P(0)=0$. For a section $V\in \mathcal C ^{\infty}(M\leftarrow \Nf)$
take the variation $f_t=f_{tV}$. The section $V$ is the variation
vector field for $f_t$ at $0$. Using now formula (\ref{Griffiths})
one sees that the linearization $L_0 \tilde P$ of $\tilde P$ at $0$
is given by the formula
\begin{equation}\label{LP}
L_0\tilde P(V)=d(\iota _W\nu ),
\end{equation} where $V=Jf_*W$ and $\nu$ is the volume form on $M$ induced by
$f$.

Since for each $V\in \mathcal C ^{\infty}(M\leftarrow \Nf)$ the
embedding $f_V$ is homotopic to $f$,  we have that $\tilde P$ has
values in $\mathcal C^{\infty}(\mathcal F_{exact}^{n})$. Moreover,
as it was observed in Section 2,  if $f_V$ is minimal affine
Lagrangian, then it is automatically special.

We shall now use the isomorphism $\rho$ given by (\ref{rho}).
   If $\gamma \in \mathcal C ^{\infty}(\mathcal F ^{n-1})$ and $V=\rho\circ\gamma$,
    then $\gamma= \iota_W\nu$, where $V=Jf_*W$.
   We now   have the mapping
$$P: \mathcal C^{\infty}({\mathcal F}^{n-1})\longrightarrow \mathcal C^{\infty}({\mathcal F}^n_{exact})$$
defined as follows:
\begin{equation}
 P(\gamma) =\tilde P (\rho \circ\gamma).
\end{equation}
The mapping $P$ can be also expressed as follows. If we set
$\Theta=({\mathcal E _f}\circ\rho)^*\Im\Omega$, then $\Theta$ is a
closed $n$-form on the  total space of $\Lambda ^{n-1}TM^*$. We have
$ P(\gamma )=\gamma ^*\Theta $ for any $(n-1)$ - form $\gamma$.
Obviously $ P(0)=0$. Moreover,  $P(\gamma)=0$ if and only if $f_V$,
where $V= \rho\circ\gamma $, is special (if $f_V$ is affine Lagrangian).

 We shall now regard $ P$ as
a differential operator. It is smooth, of order 1, non-linear and,
by (\ref{LP}), the linearization $L_0 P$ of $P$ at $0$ is given by
\begin{equation}
 L_0 P=d
\end{equation}

We shall now fix an arbitrary positive definite metric tensor field
on $M$.  The metric is  only a tool here and has no relation with
the affine geometric structure  considered in this paper. Denote by
$\delta$ the codifferential operator determined by the metric.
 Denote by ${\mathcal C}^{k,a}(\mathcal F^{r})$ the
H\"older-Banach space of all  $r$-forms on $M$ of class ${\mathcal
C}^{k,a}$, where $k\in \mathbf N$ and  $a$ is a real number from
$(0,1)$.

 We extend the action of the operators $ P, d, \delta$  to the action
 on the
  forms of class $\cka$.  The extensions will be denoted by the same letters.
In particular, after extending,  $P$ becomes a ${\mathcal
C}^{\infty}$ mapping between  Banach spaces, see \cite{Ba} p. 34,
\begin{equation}
 P: \cka (\fnj ) \longrightarrow {\mathcal
C}^{k-1,a}({\fn}_{exact})
\end{equation}
for each $k=1,2,...$. As in the proof of Theorem \ref{Mal} one sees
that there is an open neighborhood, say  $\mathcal W$, of $0$ in
$\cja (\fnj )$ such that $f_V$ are affine Lagrangian for
$V=\rho\circ\gamma$ and $\gamma \in \mathcal W$. From now on all
neighborhoods of $0$ in $\mathcal C ^{1,a}(\fnj)$ will be contained
in $\mathcal W$. Moreover, all neighborhoods will be assumed open.

Consider now  $P: \cja (\fnj ) \longrightarrow {\mathcal
C}^{0,a}({\fn}_{exact})$ as a mapping between Banach spaces. The
mapping $L_0 P$ is a surjection. Moreover $\ker L_0 P=\ker d$.
 Denote the
Banach space $\mathcal C ^{1,a}(\mathcal F_{closed}^{n-1})=\ker
d\subset\cja(\fnj)$ by $X$. The space $\delta (\cda (\fn ))$ is a
closed complement to $\ker d$. Denote this Banach space  by $Y$.
Using the implicit mapping theorem for Banach spaces one gets that
there is an open neighborhood $A$ of $0$ in $X$ and an open
neighborhood $B$ of $0$ in $Y$ and a unique smooth mapping
$G:A\longrightarrow B$ such that
$$(A+B)\cap P^{-1}(0)=\{\alpha +G(\alpha);\, \alpha \in A\}.$$

 We shall now observe that  if  $\alpha $   is of
class  ${\mathcal C} ^{k,a}$, where $k\ge 2$ or of class $\mathcal C
^{\infty}$, then so is $G(\alpha)$, for $\alpha$ from some
neighborhood of $0$ in $X$. In Riemannian geometry special submanifolds
 as minimal ones are automatically $\mathcal C ^{\infty}$
 (after possible reparametrization), but in the affine case we do not
  have such a statement and we have to prove that $\alpha$ of
   class $\mathcal C ^{\infty}$ give rise to a smooth embedding.

For an $(n-1)$--form $\gamma$  we define a differential
operator $ P_{\gamma}$ of the second order  from the vector bundle
$\Lambda ^n{TM^*}$ into itself  by the formula

\begin{equation}
P_{\gamma} (\beta)=  P(\gamma+\delta\beta)
\end{equation}
for an $n$-form $\beta$.
Since $d\beta=0$, the linearization of $ P_0$ at $0$ is the Laplace
operator.
We also have
\begin{equation}
 L_{\beta}P_{\gamma}=L_{\gamma +\delta \beta}P\circ \delta.
\end{equation}
Hence, if $\gamma$ is of class $\mathcal C ^{k,a}$ and $\beta$ is of class $\mathcal C ^{k+1,a}$, then the linear differential operator
  $L_{\beta}P_{\gamma}$ is of class  $\mathcal C ^{k-1,a}$.

We have the following smooth mapping between Banach spaces
$$\mathcal C^{1,a}(\mathcal F^{n-1})\times \mathcal C^{2,a}(\mathcal F
^{n})\ni (\gamma,\beta)\longrightarrow P_\gamma (\beta)\in \mathcal
C ^{0,a}(\mathcal F ^n)$$
and the continuous mapping
\begin{equation}
\Phi :STM^*\times \cja(\fnj)\times \mathcal C ^{2,a}(\fn)\ni (\xi, \alpha,\beta ) \longrightarrow
\det\, \sigma _{\xi} (L_{\beta} { P}_{\alpha})\in \R ,
\end{equation}
where $STM^*$ stands for the total space of the unit spheres bundle in $TM^*$ and $\sigma _\xi$ denotes
the principal symbol of a differential operator.
Since $STM^*$ is compact and $\Phi (\xi, 0,0)\ne 0$ for every $\xi \in STM^*$, we obtain the following

\begin{lemma}\label{eliptyczny}
There is a neighborhood $\mathcal U_0$ of $0$ in $\mathcal C ^{1,a}(\fnj)$ and
 a neighborhood $\mathcal V_0$ of $0$ in $\mathcal C ^{2,a}(\fn)$ such that
  for each $\gamma \in \mathcal U_0$ and $\beta \in \mathcal V_0$ the
   differential operator $L_{\beta}P_{\gamma}$ is elliptic.
\end{lemma}

From the theory of elliptic differential operators applied to
$d+\delta$ we know that the codifferential (after restricting) is a
linear homemorphism of Banach spaces
$$ \delta : \mathcal C ^{2,a}(\mathcal F ^n_{exact})\longrightarrow Y=\delta (\mathcal C ^{2,a}(\fn)).$$
Take the  neighborhood of $0$ in $X$ given by $\mathcal
U_1=G^{-1}(\delta (\mathcal V_0\cap \mathcal C ^{2,a}(\mathcal F
^n_{exact}))\cap \mathcal U_0$.  Let $\alpha \in \mathcal U_1$. Then
$G(\alpha)$ exists and there exists  $\beta\in \mathcal V_0$ such
that $G(\alpha )=\delta\beta$. Moreover $P_\alpha (\beta)=0$ and
$L_\beta P_{\alpha}$ is elliptic. Take now any $k\ge 2$ and $\alpha
\in \mathcal U_1\cap \mathcal C^{k,a}(\fnj)$. Then the differential
operator $P_\alpha$ is of class $\mathcal C ^{k-1, a}$. For
$G(\alpha)$ we have $\beta\in \mathcal V_0$ of class $\mathcal C ^2$
such that $G(\alpha) =\delta\beta$, i.e. $P_\alpha(\beta)= 0$. Hence
$\beta$ is an elliptic solution of the equation
$P_{\alpha}(\beta)=0$ and from the elliptic regularity theorem for
non-linear differential  operators we know that $\beta$ is of class
$\mathcal C ^{k+1,a}$ and consequently $G(\alpha)=\delta\beta$ is of
class $\mathcal C ^{k,a}$. Thus if $\alpha$ is of class $\mathcal C
^{\infty}$ then so is $G(\alpha)$. We have got
\begin{lemma}
  There is a neighborhood $\mathcal U_1$ of $0$ in X such that for each $k\ge 1$ we have the mapping
\begin{equation}\label{odwzorowanie_klasy_ck}
G_{| \mathcal C ^{k,a}(\fnj)}: \mathcal U_1\cap \mathcal C
^{k,a}(\fnj) \longrightarrow Y\cap \mathcal C^{k,a}(\fnj)=\delta (
\mathcal C^{k+1,a}(\mathcal F ^n)).\end{equation} Consequently we
have the mapping
\begin{equation}\label{odwzorowanie_nieskonczonosc}
G_{| \mathcal C ^{\infty}(\fnj)}: \mathcal U_1\cap\mathcal C
^{\infty}(\fnj)\longrightarrow Y\cap \mathcal C
^{\infty}(\fn)=\delta (\mathcal C ^{\infty}(\fn)).\end{equation}
\end{lemma}
We know that $G:\mathcal U_1\to Y$ is smooth between Banach spaces.
We shall now prove that the mappings (\ref{odwzorowanie_klasy_ck}),
for $k=2,...\, $, are smooth mappings between Banach spaces, when we
replace $\mathcal U_1$ by a (possibly) smaller neighborhood of $0$
in $X$. It will imply that the mapping
(\ref{odwzorowanie_nieskonczonosc}) is smooth as a mapping between
Fr\'echet spaces in a sufficiently small neighborhood of $0$ in
$\mathcal C ^{\infty}(\fnj)$.

We have the continuous mapping
\begin{equation} \mathcal C ^{1,a}(\fnj)\ni\gamma\longrightarrow (L_{\gamma}P)_{| Y}\in \mathcal L (Y,Z),\end{equation}
where $Z= \mathcal C^{0,a}(\mathcal F ^n_{exact})$ and $\mathcal L (Y,Z)$ stands
for the Banach space of continuous linear mappings from $Y$ to $Z$. We know that $(L_0P)_{| Y}=d_{|Y} :Y\to Z$
is an isomorphism (linear and topological) between Banach spaces $Y$, $Z$.  Since the set of isomorphisms is
 open in $\mathcal L (Y,Z)$, there is a neighborhood, say $\mathcal U_2$, of $0$ in $\mathcal C ^{1,a} (\fnj)$
 such that if $\gamma \in \mathcal U_2$, then $(L_{\gamma}P)_{|Y}$ is an isomorphism between $Y$ and $Z$.
  Take $\gamma \in \mathcal U_2\cap \mathcal C ^{k,a}(\fnj)$, We  have the mapping
\begin{equation}\label{mapping}
(L_{\gamma}P)_{|Y_k} : Y_k\longrightarrow \mathcal C ^{k-1,a}(\mathcal F ^n_{exact})=Z\cap \mathcal C ^{k-1,a}(\fn ),
\end{equation}
where $Y_k=Y\cap \mathcal C ^{k,a}(\fnj)=\delta (\mathcal C ^{k+1,a}(\fn))$.
As a restriction of the injection $L_{\gamma}P:Y\to Z$, it is injective.
Since $P_{| \mathcal C ^{k,a}(\fnj)}:\mathcal C ^{k,a}(\fnj)\longrightarrow \mathcal C ^{k-1, a}(\mathcal F ^n_{exact})$
 is smooth between Banach spaces, we have that
$$L_{\gamma} \left (P_{|\mathcal C ^{k,a}(\fnj)}\right): \mathcal C ^{k,a}(\fnj)\longrightarrow
 \mathcal C ^{k-1,a}(\mathcal F ^n_{exact})$$
is continuous. Hence
$$ \left(L_{\gamma} \left (P_{|\mathcal C ^{k,a}(\fnj)}\right )\right) _{|Y _k}: Y_k\longrightarrow
\mathcal C ^{k-1,a}(\mathcal F ^n_{exact})$$
is continuous.
On the other hand
$$ L_{\gamma}\left(P_{|\mathcal C ^{k,a}(\fnj)}\right) =(L_{\gamma}P)_{|\mathcal C ^{k,a}(\fnj)}.
$$
Thus the mapping given by (\ref{mapping}) is a continuous  linear
monomorphism.
 We shall now show that  it is surjective
 for  $\gamma \in \mathcal U_2\cap \mathcal U_0\cap \mathcal C ^{k,a}(\fnj)$. Let  $\mu \in \mathcal C ^{k-1,a}(\mathcal F ^n_{exact})$.
Since $\gamma\in \mathcal U_0$, by Lemma \ref{eliptyczny}, we know that $L_0P_{\gamma}$ is elliptic.  The differential operator $L_0P_{\gamma}$ is of class $\mathcal C ^{k-1,a}$.
Since $\gamma \in \mathcal U_2$, there is $\beta\in \mathcal C ^{2,a}(\fn)$ such that  $L_{\gamma} P(\delta \beta)=\mu$. From the elliptic regularity theorem we know that $\beta $ is of class $\mathcal C ^{k+1,a}$, i.e. $\delta\beta$ is of class $\mathcal C ^{k,a}$. Set $\mathcal U_3=\mathcal U_0\cap \mathcal U_2$. We have got
\begin{lemma}
There is a neighborhood $\mathcal U_3$ of $0$ in $\mathcal C ^{1,a}(\fnj)$ such that for every $\gamma\in \mathcal U_3\cap \mathcal C ^{k,a}(\fnj)$ the mapping
$$L_{\gamma}P_{|Y_k}:Y_k\longrightarrow \mathcal C ^{k-1,a}(\mathcal F ^n_{exact})$$
is an isomorphism (topological and linear).
\end{lemma}
Denote by $\tilde G:A\to\mathcal C ^{1,a}(\mathcal F ^{n-1})$ the
mapping given by $\tilde G(\alpha)=\alpha + G(\alpha)$.
 Take $\mathcal U_4=\tilde G^{-1}(\mathcal U_3)\cap \mathcal U_1\subset X$.
  Let $\alpha _0 \in \mathcal U_4\cap \mathcal C ^{k,a}(\fnj)$.
  Then $\gamma _0=\alpha _0+G (\alpha _0)$ is of class $\mathcal C ^{k,a}$
  (because $\alpha _0\in \mathcal U_1$) and $L_{\gamma_0}P_{| Y_k}:Y_k\to \mathcal C ^{k-1,a}(\mathcal F^n_{exact})$
  is an isomorphism (because $\alpha _0\in \tilde G ^{-1}(\mathcal U_3)$).
   Denote by $\tilde X _k$  the Banach space $\ker L_{\gamma _0}P$.
    We have $P(\gamma _0)=0$ and $\tilde X _k\oplus Y_k=\mathcal C ^{k,a}(\mathcal F ^{n-1})$.
     We want to prove that $G$ is smooth  around $\alpha _0$ in the sense of the Banach spaces theory.
      Denote by $\tilde \pi :\mathcal C ^{k,a}(\fnj)=\tilde X _k\oplus Y_k\longrightarrow \tilde X_k$
      the canonical projection. It is a smooth mapping between Banach spaces.
      Set $\tilde \alpha _0=\tilde \pi(\alpha _0)$.
      From the implicit mapping theorem we know that there is a neighborhood $\tilde U $ of $\tilde\alpha _0$
      in $\tilde X _k$ and a smooth mapping $F$ defined on $\tilde U$ such that
       $ \{\tilde\alpha + F(\tilde\alpha); \tilde\alpha \in \tilde U\}\subset P^{-1}(0)$.
 In a neighborhood of $\alpha _0$ we have
$$G(\alpha)= \tilde\pi (\alpha) +F(\tilde \pi (\alpha))-\alpha.$$
Hence in a neighborhood of $\alpha _0$ the mapping $G$ is smooth and
consequently it is smooth in $\mathcal U_4\cap \mathcal C
^{k,a}(\fnj)$. It follows that $G_{|\mathcal U_4\cap \mathcal C
^{\infty}(\fnj)}$ is smooth in the sense of the theory of Fr\'echet
 spaces.
   The projections in the
Hodge decomposition $\mathcal C ^{\infty}(\mathcal F^{n-1})=
\mathcal C^{\infty}(\mathcal F^{n-1}_{closed})\oplus  \delta
(\mathcal C ^{\infty}(\mathcal F ^n))$ are smooth mappings of
Fr\'echet spaces.
Denote by $p: \mathcal C ^{\infty}(\mathcal F^{n-1})\to \mathcal
C^{\infty}(\mathcal
 F^{n-1}_{closed})$
the projection. Set $\mathfrak U=(\mathcal U_4 \cap \mathcal C
^{\infty}(\mathcal F^{n-1}))\oplus (\delta(\mathcal C^{\infty
}\mathcal (F^{n}))$. Consider the mapping $\phi: \mathfrak U\to
\mathfrak U$ defined as $\phi (z)=z-(G\circ p)(z)$. It is a
bijection and its converse  is given by  $\phi ^{-1}(z)= z+(G\circ
p)(z)$. Both mappings $\phi$ and $\phi ^{-1}$ are smooth  in the
sense of the theory of Fr\'echet vector spaces. We  now compose the
chart obtained in the proof of Theorem \ref{Mal}
  with $\phi$. Since
$$\phi ( \{\gamma =\alpha +G(\alpha); \
\alpha\in \mathcal U_4\cap \mathcal C ^{\infty}(\mathcal F ^{n-1})
\})=\mathcal U_4 \cap \mathcal C ^{\infty}(\mathcal F ^{n-1})$$
is an open subset of the closed subspace $\mathcal C
^{\infty}(\mathcal F^{\infty}_{closed})$ of $\mathcal C
^{\infty}(\mathcal F ^{\infty})$, we have that the set $\mathcal M
maL$ is a submanifold of $\mathcal M aL$.
 The proof is completed.

\bigskip
\begin{remark}
{\rm We shall now observe that there exist minimal affine Lagrangian submanifolds
which are not smooth. We refer to Section 3 for notation.
 Having a smooth special affine Lagrangian embedding $f:M\to N$ we have the mapping
$$ \Phi _k:\mathcal C ^k_{emb}(M, \mathcal T _f)\ni h\longrightarrow \Pi _f
\circ \mathcal E ^{-1}_f\circ h\in \mathcal C ^k(M,M).$$
The set $Diff ^k(M)$ is open in $\mathcal C ^k(M,M)$ and the set
$$\mathcal U _{k,f}=\{ h\in \mathcal C ^k_{emb}(M,\mathcal T _f) ;\
\exists V\in \mathcal C ^k(M\leftarrow \mathcal N _f): \mathcal E _f\circ V=h\}$$
can be regarded as an open neighbourghood of $[f]$ in
$\mathcal M^ {k}=\mathcal C ^k(M, \mathcal T _f)_{/ Diff ^k(M)}$.
We have the bijection
\begin{equation}\label{bijection}
\mathcal C ^k(M,\leftarrow \mathcal N_f) \ni V\longrightarrow \mathcal E _f\circ V\in \mathcal U  _{k,f}.
\end{equation}
In order to study minimal affine Lagrangian submanifolds of complex
equiaffine spaces like in \cite{O_1}, \cite{O} or in Section 1 of
this paper it suffices that the immersions or embeddings under
consideration are of class $\mathcal C ^2$. Also in the  proof of
Theorem \ref{main-theorem} the class $\mathcal C ^2$ is sufficient,
that is, if $\alpha$ is of class $\mathcal C ^k$, where $k\ge 2$,
then $G(\alpha)$ is of class $\mathcal C ^k$. Since $\mathcal C
^k(\mathcal F ^{n-1}_{closed})\ne\mathcal C ^{\infty}(\mathcal F
^{n-1}_{closed})$, it is clear by the proof of Theorem
\ref{main-theorem} that there exist non-smooth minimal affine
Lagrangian embeddings of class $\mathcal C ^k$, for $k\ge 2$. }

\end{remark}
\medskip

\begin{example}{\rm Let $M$ be an $n$-dimensional real manifold equipped with a torsion-free linear connection
$\nabla$.  The tangent bundle to the tangent bundle $TTM$ admits a
decomposition into a direct sum of the vertical bundle (tangent to
the fibers of $TM$) and the  horizontal bundle (depending on the
connection). The vertical lift of $X\in T_xM$ to $T_ZTM$ for $Z\in
T_xM$ will be denoted by $X^v_Z$. Analogously the horizontal lift
will be denoted by $X^h_Z$. The following  formulas for the lifts of
vector fields $X, Y\in \mathcal X (M)$ are known, see \cite{D},

\begin{equation}\label{nawiasy}
\begin{array}{rcl}
&& [X^v,Y^v]=0,\\
&&[X^h, Y^v]=(\nabla _XY)^v,\\
&& [X^h,Y^h]_Z=-(R(X,Y)Z)^v_Z +[X,Y]^h_Z,
\end{array}
\end{equation}
where $R$ denotes the curvature tensor of $\nabla$.

The total space  $TM$ has an almost complex structure $J$ determined
by $\nabla$. Namely
\begin{equation}
JX^h=X^v,\ \ \ \ \ \ JX^v=-X^h.
\end{equation}
 From
(\ref{nawiasy}) it follows that the almost complex structure is
integrable if and only if the connection $\nabla$ is flat.

Assume that $\nu$ is  a volume form on $M$ such that $\nabla \nu
=0$. In other words, the pair $\nabla$, $\nu$ is an equiaffine
structure on $M$. We define a complex volume form $\Omega$ on $TM$
by the formula
\begin{equation}\label{omega}
\Omega (X_1^h,...,X_n^h)=\nu (X_1,...,X_n).
\end{equation}
By using (\ref{nawiasy}) one sees that $d\Omega=0$ if and only if
$\nabla$ is flat.

From now on we assume that $\nabla$ is flat and $\nabla\nu =0$. A
manifold with such a structure is usally called an  affine manifold
with parallel volume. Take the zero-section of $TM$.
 The horizontal space at $0_x$ is equal to $T_xM$
(independently of a connection $\nabla$). Hence the zero-section
treated as a mapping $ 0:M\to TM$ is an affine Lagrangian embedding.
By (\ref{omega}) it is special (also independently of a given
connection). We have
\begin{proposition}
Each affine manifold with parallel volume admits a special affine
Lagrangian embedding into a complex space with closed  complex
volume form.
\end{proposition}
From the main theorem of this paper we know that if $M$ is
additionally compact, then such embeddings are plentiful.

If $\nabla$ is flat and $\nabla\nu =0$ then, in fact,  the total
space of the tangent bundle $TM$ has a structure of a complex
equiaffine manifold.}
\end{example}

\begin{remark}
{\rm Assume now that $N$ is a 4-dimensional almost K\"ahler manifold
with symplectic form $\kappa$. Let $M$ be a connected compact
orientable $2$-dimensional manifold and $f:M\to N$ be a Lagrangian
embedding (in the metric sense). We now have the canonical
(depending only on the metric) isomorphism, say $\mathfrak{b}$,
between vector fields and $1$-forms on $M$. By Theorem  \ref{Mal} we
have the manifold $\mathcal M aL$ modeled on the Fr\'echet space
$\mathcal C ^{\infty}(\mathcal F ^1)$.

Similarly as in the proof of Theorem \ref{main-theorem} we define
the mapping

\begin{equation}
\tilde P : \mathcal C ^{\infty}(M\leftarrow \mathcal N) \ni V\to
f^*_V\kappa \in\mathcal C ^{\infty}(\mathcal F^2).
\end{equation}
Since $f^*\kappa =0$ and the normal bundle is star-shaped, the
mapping $\tilde P$ has values in $\mathcal C ^{\infty}(\mathcal
F^2_{exact})$. By composing the mapping with the isomorphism between
the normal bundle and the tangent bundle and the isomorphism
$\mathfrak{b}$ we obtain the mapping $$ P:\mathcal C
^{\infty}(\mathcal F^1)\to \mathcal C ^{\infty}(\mathcal F ^2)$$
whose linearization at $0$ is equal to the exterior differential
operator $d$. Now we can argue as in the proof of Theorem
\ref{main-theorem} and we get

\begin{thm}
Let $N$ be a 4-dimensional almost K\"ahler manifold and $f:M\to N$
be a Lagrangian embedding of a connected compact orientable
2-dimensional manifold. Then the set $$\mathcal M L= \{ [f]\in
\mathcal MaL;\ f\ is \ Lagrangian\}$$ is an infinite dimensional
Fr\'echet manifold modeled on the Fr\'echet vector space $\mathcal
C^{\infty}(\mathcal F ^1_{closed})$. It is a submanifold of
$\mathcal M aL$.
\end{thm}

}
\end{remark}

\end{document}